\documentclass[11pt,amsfonts]{article}
\usepackage{latexsym}
\usepackage{amssymb}
\usepackage{amsmath}
\usepackage{layout}

\newtheorem{lemma}{Lemma}

\newtheorem{theorem}{Theorem}

\def\real{{\mathord{{\rm I\kern-2.8pt R}}}}        % Fake blackboard bold R.

\def\inte{{\mathord{{\rm I\kern-2.8pt N}}}}

\def\sZZ{{\rm Z\kern-2.8ptem{}Z}}

\def\z{{\mathchoice
  {\sZZ}
  {\sZZ}
  {\rm Z\kern-0.30em{}Z}
  {\rm Z\kern-0.25em{}Z} }}
\def\sQQ{{\kern 0.27em \vrule height1.45ex width0.03em depth0em
          \kern-0.30em \rm Q}}
\def\qu{{\mathchoice
    {\sQQ}
    {\sQQ}
  {\kern 0.225em \vrule height1.05ex width0.025em depth0em \kern-0.25em
\rm Q}
  {\kern 0.180em \vrule height0.78ex width0.020em depth0em \kern-0.20em
\rm Q}
        }}
\def\sCC{{\kern 0.27em \vrule height1.45ex width0.03em depth0em
          \kern-0.30em \rm C}}
\def\complex{{\mathchoice
    {\sCC}
    {\sCC}
  {\kern 0.225em \vrule height1.05ex width0.025em depth0em \kern-0.25em
\rm C}
  {\kern 0.180em \vrule height0.78ex width0.020em depth0em \kern-0.20em
\rm C}
        }}

\newcommand{\R}{\mathbb{R}}

\newcommand{\ba}{\begin{array}}
\newcommand{\ea}{\end{array}}
\newcommand{\be}{\begin{equation}}
\newcommand{\ee}{\end{equation}}
\newcommand{\bea}{\begin{eqnarray}}
\newcommand{\eea}{\end{eqnarray}}
\newcommand{\beaa}{\begin{eqnarray*}}
\newcommand{\eeaa}{\end{eqnarray*}}

\def\z{\zeta}
\def\k{\kappa}

\font\tenmath=msbm10 \font\sevenmath=msbm7 \font\fivemath=msbm5
\newfam\mathfam \textfont\mathfam=\tenmath
\scriptfont\mathfam=\sevenmath \scriptscriptfont\mathfam=\fivemath
\def\math{\fam\mathfam}

\def \={{\buildrel {\rm (law)} \over =}}

\def \R{{\math R}}

\def\qed{ \hfill \vrule width.25cm height.25cm depth0cm\smallskip}

\newcommand{\basa}{\begin{assumption}}
\newcommand{\easa}{\end{assumption}}

\newcommand{\bas}{\begin{assum}}
\newcommand{\eas}{\end{assum}}

\newcommand{\ignore}[1]{}
\begin{document}

\renewcommand{\thefootnote}{\fnsymbol{footnote}}

\title{Multidimensional bifractional Brownian motion: It\^o and Tanaka
formulas}
\author{Khalifa Es-sebaiy  $^{1}\quad$Ciprian A. Tudor
$^{2}\vspace*{0.1in}$\\$^{1}$
  Department of Mathematics, Faculty
of Sciences Semlalia,\\ Cadi Ayyad University 2390 Marrakesh,
Morocco.
\\k.essebaiy@ucam.ac.ma\vspace*{0.1in}\\$^{2}$SAMOS/MATISSE,
Centre d'Economie de La Sorbonne,\\ Universit\'e de
Panth\'eon-Sorbonne Paris 1,\\90, rue de Tolbiac, 75634 Paris Cedex
13, France.\\tudor@univ-paris1.fr\vspace*{0.1in}}  \maketitle

\maketitle

 \begin{abstract}
Using the Malliavin calculus with respect to Gaussian processes and
the multiple stochastic integrals we derive It\^o's and Tanaka's
formulas for the $d$-dimensional bifractional Brownian motion.
\end{abstract}

 \vskip0.5cm

{\bf  2000 AMS Classification Numbers: 60G12, 60G15,  60H05, 60H07.
}

 \vskip0.3cm

{\bf Key words:} Malliavin calculus, bifractional Brownian motion,
multiple stochastic integrals, Watanabe distribution, Tanaka
formula.

\vskip0.3cm

\section{Introduction}

The stochastic calculus with respect to the fractional Brownian
motion (fBm) has now a long enough history. Since the nineties, many
authors used different approaches to develop a stochastic
integration theory with respect to this process. We refer, among of
course many others, to \cite{AMN}, \cite{DU}, \cite{RV} or
\cite{DHD}. The reason for this tremendous interest in the
stochastic analysis of the fBm comes from its large amount of
applications in practical phenomena such as telecommunications,
hydrology or economics.

Nevertheless, even fBm  has its limits in modeling  certain phenomena. Therefore, several authors introduced recently some
generalizations of the fBm which are supposed to fit better in concrete situations. For example, we mention the
multifractional Brownian motion (see e.g. \cite{AY}), the subfractional Brownian motion (see e.g. \cite{Gor2}) or the
multiscale fractional Brownian motion (see \cite{BB}).

Here our main interest consists in the study of the {\em
bifractional Brownian motion (bifBm)}. The bifBm has been introduced
by Houdr\'e and Villa in \cite{HV} and a stochastic analysis for it
can be found in \cite{russo}. Other papers  treated different
aspects of this stochastic process, like sample pats properties,
extension of the parameters or statistical applications (see
\cite{Gor}, \cite{Beg}, \cite{Xiao} or  \cite{chine}. Recall that
the bifBm $B^{H,K}$ is a centered Gaussian process, starting from
zero, with covariance function
\begin{equation}
\label{cov}R^{H,K}(t,s) := R(t,s)= \frac{1}{2^{K}}\left( \left(
t^{2H}+s^{2H}\right)^{K} -\vert t-s \vert ^{2HK}\right)
\end{equation}
where the parameters $H,K$ are such that $H\in (0,1)$ and $K\in
(0,1]$. In the case $K=1$ we retrieve  the fractional Brownian
motion while  the case $K=1$ and $H=\frac{1}{2}$ corresponds to the
standard Brownian motion.

The process $B^{H,K}$ is $HK$-selfsimilar but it has no stationary
increments. It has H\"older continuous paths of order $\delta <HK$
and its paths are not differentiable. An interesting property of it
is the fact that its quadratic variation in the case $2HK=1$ is
similar to that of the standard Brownian motion, i.e. $[B^{H,K}]_{t}
=cst. \times t $ and therefore especially this case ($2HK=1$) is
very interesting from the stochastic calculus point of view.

In this paper, our purpose is to study multidimensional bifractional Brownian motion and to prove It\^o and Tanaka formulas.
We start with the one dimensional bifBm and  we first derive an It\^o and an Tanaka formula for it when $2HK\geq 1$. We
mention that the It\^o formula has been already proved by \cite{KRT} but here we propose an alternative proof based on the
Taylor expansion which appears to be also useful in the multidimensional settings. The Tanaka formula is obtained from the
It\^o formula by a limit argument and it involves the so-called {\em weighted local time} extending the result in \cite{CNT}.
In the multidimensional case we first derive an It\^o formula for $2HK>1$ and we extend it to Tanaka by following an idea by
Uemura \cite{uemura1}, \cite{uemura}; that is, since $\vert x\vert$ is  twice the kernel of the one-dimensional Newtonian
potential, i.e. $\frac{1}{2}\Delta \vert x\vert $ is equal to the delta Dirac function $\delta (x)$, we will chose the
function $U(z), z \in \mathbb{R}^{d}$ which is twice of the kernel of $d$-dimensional Newtonian (or logarithmic if $d=2$)
potential to replace $\vert x\vert$ in the $d$-dimensional case. See the last section for the definition of the function $U$.
Our method is based on the Wiener-It\^o chaotic expansion into multiple stochastic integrals following ideas from \cite{Im} or
\cite{ELSTV}. The multidimensional Tanaka formula also involves a generalized local time. We note that the terms appearing in
our Tanaka formula when $d\geq 2$ are not random variables and they are understood as distributions in the Watanabe spaces.

\section{Preliminaries: Deterministic spaces associated and
Malliavin calculus} Let $\left( B^{H,K}_{t}, t\in [0,T]\right) $ be
a bifractional Brownian motion on the probability space $(\Omega,
{\cal{F}}, P)$.

Being a Gaussian process, it is possible to construct a stochastic
calculus of variations with respect to $B^{H,K}$. We refer to
\cite{AMN}, \cite{N} for a complete description of stochastic
calculus with respect to Gaussian processes. Here we recall only the
basic elements of this theory.

The basic ingredient is the canonical Hilbert space ${\cal{H}}$
associated to the bifractional Brownian motion (bifBm).  This space
is defined as the completion of the linear space ${\cal{E}}$
generated by the indicator functions ${1_{[0,t]}, t\in [0,T]}$  with
respect  to the inner product
\begin{equation*}
\langle 1_{[0,t]}, 1_{[0,s]}\rangle _{{\cal{H}}} =R(t,s).
\end{equation*}
The application $\varphi \in {\cal{E}}\to B(\varphi )$ is an
isometry from ${\cal{E}}$ to  the Gaussian space generated by
$B^{H,K}$ and it can be extended to ${\cal{H}}$.

Let us denote by ${\cal{S}}$ the set of smooth functionals of the
form
\begin{equation*}
F=f(B(\varphi _{1}), \ldots, B(\varphi _{n}))
\end{equation*}
where $f\in C ^{\infty }_{b}(\mathbb{R}^{n})$  and $\varphi _{i} \in
{\cal{H}}$. The Malliavin derivative of a functional $F$ as above is
given by
\begin{equation*}
D^{B^{H,K}}F= \sum _{i=1}^{n} \frac{\partial f}{\partial
x_{i}}(B(\varphi _{1}), \ldots, B(\varphi _{n})) \varphi _{i}
\end{equation*}
and this operator can be extended to the closure $\mathbb{D} ^{m,2}$
($m\geq 1$) of ${\cal{S}}$ with respect to the norm
\begin{equation*}
\Vert F\Vert ^{2}_{m,2} := E\left| F \right| ^{2} + E\Vert
D^{B^{H,K}}F\Vert ^{2}_{{\cal{H}}}+\ldots + E\Vert D^{B^{H,K},m}
F\Vert ^{2} _{{\cal{H}}^{ \hat{\otimes } m}}
\end{equation*}
where ${\cal{H}}^{ \hat{\otimes } m}$ denotes the $m$ fold symmetric
tensor product and the $m$th derivative $D^{B^{H,K},m}$ is defined
by iteration.

The divergence integral $\delta ^{B^{H,K}}$ is the adjoint operator of $D^{B^{H,K}}$. Concretely, a random variable $u\in
L^{2}(\Omega ; {\cal{H}}) $ belongs to the domain of the divergence operator ($Dom(\delta ^{B^{H,K}})$) if
$$E\left|  \langle D^{B^{H,K}}F, u\rangle _{{\cal{H}}} \right|\leq c \Vert F\Vert
_{L^{2}(\Omega )}  $$ for every $F\in {\cal{S}}$. In this case $\delta ^{B^{H,K}} (u)$ is given by the duality relationship
\begin{equation*}
E(F\delta ^{B^{H,K}} (u)) = E\langle D^{B^{H,K}}F, u\rangle _{{\cal{H}}}
\end{equation*}
for any $F\in \mathbb{D}^{1,2}$. It holds that
\begin{equation}
\label{delta 2} E\delta ^{B^{H,K}} (u) ^{2} = E\Vert u \Vert ^{2} _{ {\cal{H}}} + E\langle D^{B^{H,K}}u ,
(D^{B^{H,K}}u)^{\ast} \rangle _{{\cal{H}}\otimes {\cal{H}}}
\end{equation}
where $(D^{B^{H,K}}u)^{\ast}$ is the adjoint of $D^{B^{H,K}}u$ in
the Hilbert space ${\cal{H}}\otimes {\cal{H}}$.

\vskip0.3cm

Sometimes working with the space ${\cal{H}}$ is not convenient;
once, because this space may contain also distributions (as, e.g. in
the case $K=1$, see \cite{PT}) and twice, because the norm in this
space is not always tractable. We will use the subspace $\left|
{\cal{H}}\right| $ of ${\cal{H}}$ which is defined   as the set of
measurable function $f$ on $[0,T]$ with
\begin{equation}
\label{barH} \Vert f\Vert ^{2} _{  \left| {\cal{H}}\right|   } :=
\int_{0}^{T} \int_{0}^{T} \left| f(u) \right| \left| f(v) \right|
\frac{\partial ^{2}R}{\partial u\partial v} (u,v) dudv <\infty.
\end{equation}
It follows actually from \cite{KRT} that the space $\left|
{\cal{H}}\right|$ is a Banach space for the norm $\Vert \cdot \Vert
 _{  \left| {\cal{H}}\right|   }$ and it is included in
${\cal{H}}$. In fact,
\begin{equation}
\label{inclu} L^{2}([0,T])\subset L^{\frac{1}{HK}}([0,T]) \subset
\left| {\cal{H}}\right| \subset {\cal{H}}.
\end{equation}
and
\begin{equation}
\label{ineg1} E\delta ^{B^{H,K}} (u)^{2} \leq E\Vert u \Vert  ^{2} _{\left| {\cal{H}}\right| } + E\Vert D^{B^{H,K}}u\Vert ^{2}
_{ \left| {\cal{H}}\right|\otimes \left| {\cal{H}}\right| }
\end{equation}
where, if $\varphi :[0,T]^{2}\to \mathbb{R}$
\begin{equation}\label{HH}
 \Vert \varphi \Vert ^{2} _{ \left|
{\cal{H}}\right|\otimes \left| {\cal{H}}\right| } = \int
_{[0,T]^{4}} \left|  \varphi (u,v)\right| \left| \varphi (u', v')
\right| \frac{\partial ^{2} R}{\partial u\partial
u'}(u,u')\frac{\partial ^{2} R}{\partial v\partial
v'}(v,v')dudvdu'dv'.
\end{equation}
We will use the following formulas of the Malliavin calculus: the
integration by parts
\begin{equation}
\label{ip} F\delta ^{B^{H,K}} (u)= \delta ^{B^{H,K}} (Fu) + \langle D^{B^{H,K}}F, u\rangle _{ {\cal{H}}}
\end{equation}
for any $u\in Dom(\delta ^{B^{H,K}})$, $ F\in \mathbb{D}^{1,2}$ such that $Fu\in L^{2}(\Omega; {\cal{H}})$; and the chain rule
\begin{equation}
\label{chain} D^{B^{H,K}} \varphi(F)= \sum _{i=1}^{n} \partial
_{i}\varphi (F) D^{B^{H,K}}F^{i}
\end{equation}
if $\varphi : \mathbb{R} ^{n}\to \mathbb{R}$ is continuously
differentiable with bounded partial derivatives and $F=(F^{1},
\ldots , F^{m})$ is a random vector with components in
$\mathbb{D}^{1,2}$.

\vskip0.3cm

 By the duality between $D^{B^{H,K}}$ and $\delta ^{B^{H,K}}$ we obtain the
following result for the convergence of divergence integrals: if $u_{n}\in Dom (\delta ^{B^{H,K}})$ for every $n$,
$u_{n}\underset{n\to \infty}{\to}u$ in $L^{2}(\Omega ; {\cal{H}})$ and $\delta ^{B^{H,K}} (u_{n} ) \underset{n\to \infty}{\to}
G$ in $L^{2}(\Omega)$ then
\begin{equation}
\label{basic} u\in Dom(\delta ^{B^{H,K}}) \mbox{ and } \delta ^{B^{H,K}} (u)= G.
\end{equation}

It is also possible to introduce multiple integrals $I_{n}(f_{n})$,
$f\in {\cal{H}}^{\otimes n}$ with respect to $B^{H,K}$. Let
\begin{equation}\label{sum1} F=\sum _{n\geq 0} I_{n}(f_{n})
\end{equation}
where for every $n\geq 0$, $f_{n} \in {\cal{H}}^{\otimes n}$ are
symmetric functions. Let $L$ be the Ornstein-Uhlenbeck operator
\begin{equation*}
LF=-\sum_{n\geq 0} nI_{n}(f_{n})
\end{equation*}
if $F$ is given by (\ref{sum1}).

For $p>1$ and $\alpha \in \mathbb{R}$ we introduce the
Sobolev-Watanabe space $\mathbb{D}^{\alpha ,p }$  as the closure of
the set of polynomial random variables with respect to the norm
\begin{equation*}
\Vert F\Vert _{\alpha , p} =\Vert (I -L) ^{\frac{\alpha }{2}} \Vert
_{L^{p} (\Omega )}
\end{equation*}
where $I$ represents the identity. In this way, a random variable
$F$ as in (\ref{sum1}) belongs to $\mathbb{D}^{\alpha , 2}$ if and
only if
\begin{equation*}
\sum_{n\geq 0} (1+n) ^{\alpha } \Vert I_{n}(f_{n}) \Vert
_{L^{2}(\Omega)} ^{2}<\infty.
\end{equation*}
Note that  the Malliavin  derivative operator  acts on multiple
integral as follows
\begin{equation*}
D^{B^{H,K}}_{t} F=\sum_{n=1}^{\infty}nI_{n-1}(f_{n}(\cdot , t)),
\hskip0.5cm t\in [0,T].
\end{equation*}
The operator $D^{B^{H,K}}$ is continuous from $\mathbb{D} ^{\alpha -1 , p} $ into $\mathbb{D} ^{\alpha , p} \left(
{\cal{H}}\right).$ The adjoint of $D^{B^{H,K}}$ is denoted by $\delta ^{B^{H,K}} $ and is called the divergence (or Skorohod)
integral. It is a continuous operator from $\mathbb{D} ^{\alpha, p } \left( {\cal{H}}\right)$ into $\mathbb{D} ^{\alpha , p}$.
For adapted integrands, the divergence integral coincides to the classical It\^o integral. We will use the notation
\begin{equation*}
\delta ^{B^{H,K}} (u) =\int_{0}^{T} u_{s} \delta  B_{s}.
\end{equation*}
\vskip0.5cm

 Recall
that if $u$ is a stochastic process having the chaotic decomposition
\begin{equation*}
u_{s}=\sum _{n\geq 0} I_{n}(f_{n}(\cdot ,s))
\end{equation*}
where $f_{n}(\cdot, s)\in {\cal{H}}^{\otimes n}$ for every $s$, and
it is symmetric in the first $n$ variables, then its Skorohod
integral is given by
\begin{equation*}
\int_{0}^{T} u_{s}dB_{s}= \sum _{n\geq 0}I_{n+1}(\tilde{f}_{n})
\end{equation*}
where $\tilde{f}_{n}$ denotes the symmetrization of $f_{n}$ with
respect to all $n+1$ variables.

\section{Tanaka formula for unidimensional bifractional Brownian
motion}

This paragraph is consecrated to the proof of It\^o formula and
Tanaka formula for the one-dimensional bifractional Brownian motion
with $2HK\geq 1$. Note that the It\^o formula has been already
proved in \cite{KRT}; here we propose a different approach based on
the Taylor expansion which be also used in the multidimensional
settings.

\vskip0.3cm

We start by the following technical lemma.

\begin{lemma}\label{lem1} Let us consider the following function on
$[1,\infty)$
\begin{eqnarray*}h(y)=y^{2HK}+(y-1)^{2HK}-\frac{2}{2^K}\left(y^{2H}+(y-1)^{2H}\right)^K.
\end{eqnarray*} where $H\in(0,1)$ and $K\in(0,1)$.
Then, \begin{eqnarray}h(y)\  \mbox{ converges to } 0 \mbox{ as } y
\mbox{ goes to } \infty.\label{function h}\end{eqnarray}Moreover if
$2HK=1$ we obtain that
\begin{eqnarray}
\lim_{y\rightarrow+\infty}yh(y)=\frac{1}{4}(1-2H).\label{2HK=1}
\end{eqnarray}
\end{lemma}
\begin{proof}  Let
$\displaystyle y=\frac{1}{\varepsilon} $, hence
\begin{eqnarray*}h(y)=h\left(\frac{1}{\varepsilon}\right)
=\frac{1}{\varepsilon^{2HK}}\left[1+(1-\varepsilon)^{2HK}-\frac{2}{2^K}\left(1+(1-\varepsilon)^{2H}\right)^K\right]
\end{eqnarray*}
Using Taylor's expansion, as $\varepsilon$ close to 0, we obtain
\begin{eqnarray}h\left(\frac{1}{\varepsilon}\right)
=\frac{1}{\varepsilon^{2HK}}\left(H^2K(K-1)\varepsilon^{2}+o(\varepsilon^{2})\right)
\label{taylor}
\end{eqnarray}
Thus
\begin{eqnarray*}
\lim_{y\rightarrow+\infty}h(y)=\lim_{\varepsilon\rightarrow0}h(1/\varepsilon)=0.
\end{eqnarray*}
For the case $2HK=1$ we replace in (\ref{taylor}), we have
\begin{eqnarray*}\frac{1}{\varepsilon}h\left(\frac{1}{\varepsilon}\right)
=\frac{1}{4}(1-2H)+\frac{1}{\varepsilon^{2}}o(\varepsilon^{2})
\end{eqnarray*} Thus (\ref{2HK=1}) is satisfied. Which
completes the proof. \qed \end{proof}

\vskip0.5cm

\begin{theorem}\label{ito unidim}Let $f$ be a
function of class $\mathrm{C}^2$ on $\mathbb{R}$. Suppose that
$2HK\geq 1$, then
\begin{eqnarray}
f\left(B^{H,K}_{t}\right)=f\left(0\right)+\int_{0}^{t}f'\left(B^{H,K}_{s}\right)\delta  B^{H,K}_{s}+HK
\int_{0}^{t}f''\left(B^{H,K}_{s}\right)s^{2HK-1}ds.\label{ito}
\end{eqnarray}
\end{theorem}
\begin{proof}
 We first prove the case $2HK>1$. Let us fix $t>0$ and let be
$\pi: =\{t_j=\frac{jt}{n}; j=0,...,n\}$ a partition of $[0,t]$. By
the localization argument and the fact that the process $B^{H,K}$ is
continuous, We can assume that f has compact support, and so $f$,
$f'$ and $f''$ are bounded. Using Taylor expansion, we have
\begin{eqnarray}f\left(B^{H,K}_t\right)&=&f(0)+\sum_{j=1}^{n}f'\left(B^{H,K}_{t_{j-1}}\right)
\left(B^{H,K}_{t_{j}}-B^{H,K}_{t_{j-1}}\right)+\frac{1}{2}\sum_{j=1}^{n}f''\left(\overline{B}^{H,K}_{j}\right)
\left(B^{H,K}_{t_{j}}-B^{H,K}_{t_{j-1}}\right)^2\nonumber\\&:=&f(0)+
I^n + J^n.\label{I+J}
\end{eqnarray}
where
$\overline{B}^{H,K}_{j}=B^{H,K}_{t_{j-1}}+\theta_j\left(B^{H,K}_{t_{j}}-B^{H,K}_{t_{j-1}}\right)$
with $\theta_j$ is a r.v in $(0,1)$. \\
Since $B^{H,K}$ is a  quasi-helix (see \cite{russo}), we can bound
the term $J^{n}$ as follows:
\begin{eqnarray*}E|J^n|^2&\leq&
nC/4\sum_{j=1}^{n}E\left(B^{H,K}_{t_{j}}-B^{H,K}_{t_{j-1}}\right)^4\\
&\leq&  2^{-2K}nC\sum_{j=1}^{n}|t_j-t_{j-1}|^{4HK} \\&\leq&
2^{-2K}nC t^{4HK}\sup_{j=1}^{n}|t_j-t_{j-1}|^{4HK-1}=2^{-2K}C
\frac{t^{4HK}}{n^{4HK-2}}\underset{n\rightarrow\infty}{\longrightarrow}0.
\end{eqnarray*}where $C$ a constant depends de $f''$.
Then
\begin{eqnarray}J^n\underset{n\rightarrow\infty}{\longrightarrow}0
\mbox{ in } L^2(\Omega).\label{J}
\end{eqnarray} On the other hand, we apply (\ref{ip}) and we get
\begin{eqnarray*}I^n&=&\sum_{j=1}^{n}f'\left(B^{H,K}_{t_{j-1}}\right)\left(\delta^{B^{H,K}}(1_{(t_{j-1},t_j]})\right)\\
&=&\delta^{B^{H,K}}\left(\sum_{j=1}^{n}f'\left(B^{H,K}_{t_{j-1}}\right)1_{(t_{j-1},t_j]}(.)\right)+
\sum_{j=1}^{n}f''\left(B^{H,K}_{t_{j-1}}\right)\langle1_{(0,t_{j-1}]},1_{(t_{j-1},t_j]}\rangle_{\mathcal{H}}\\&=&
I_1^n+I_2^n.
\end{eqnarray*}
Next
\begin{eqnarray*}I^n_2&=&\sum_{j=1}^{n}f''\left(B^{H,K}_{t_{j-1}}\right)\left(R(t_{j-1},t_{j})-R(t_{j-1},t_{j-1})\right)\\
&=&\sum_{j=1}^{n}\left[f''\left(B^{H,K}_{t_{j-1}}\right)\left(\frac{1}{2^K}\left((t_j^{2H}+t_{j-1}^{2H})^K
-(t_j-t_{j-1})^{2HK}\right)-t_{j-1}^{2HK}\right)\right]
\end{eqnarray*}
We denote by
\begin{eqnarray*}A_t:=HK\int_{0}^{t}s^{2HK-1}ds=\frac{1}{2}t^{2HK}
\end{eqnarray*}
To prove that $I^n_2$ converges to
$HK\int_{0}^{t}f''\left(B^{H,K}_s\right)s^{2HK-1}ds$ in
$L^2(\Omega)$ as $n\rightarrow\infty$, it suffices to show that
\begin{eqnarray*}C_n:=\left(E\left|I^n_2-
\sum_{j=1}^{n}f''\left(B^{H,K}_{t_{j-1}}\right)\left(A_{t_{j}}-A_{t_{j-1}}\right)\right|^2\right)^{1/2}
\underset{n\rightarrow\infty}{\longrightarrow}0.
\end{eqnarray*}

 By Minkowski   inequality, we have
\begin{eqnarray*}C_n&\leq&C\sum_{j=1}^{n}
 \left|\left(\frac{1}{2^K}((t_j^{2H}+t_{j-1}^{2H})^K
-(t_j-t_{j-1})^{2HK})-t_{j-1}^{2HK}\right)-\frac{1}{2}(t_{j}^{2HK}-t_{j-1}^{2HK})\right|
\\&\leq&Ct^{2HK}/2\left[\frac{1}{n^{2HK}}\sum_{j=1}^{n}\left|h(j)\right|
+\frac{2}{2^{K}}\frac{1}{n^{2HK-1}}\right]\\&=&Ct^{2HK}/2\left[C_n^1+C_n^2\right]
\end{eqnarray*}
where $C$ is a  generic constant.\\ Since $2HK>1$ then
$C_n^2:=\frac{2}{2^{K}}\frac{1}{n^{2HK-1}}\underset{n\rightarrow\infty}{\longrightarrow}0.$
According to (\ref{function h}), we obtain
\begin{eqnarray*}C_n^1:=
\frac{1}{n^{2HK}}\sum_{j=1}^{n}h(j)\leq \frac{C}{n^{2HK-1}}
\underset{n\rightarrow\infty}{\longrightarrow}0.
\end{eqnarray*}
Thus \begin{eqnarray*}I^n_2
\underset{n\rightarrow\infty}{\longrightarrow}HK\int_{0}^{t}f''\left(B^{H,K}_s\right)s^{2HK-1}ds
\mbox{ in }
 L^2(\Omega) \end{eqnarray*}
 We show now that  $$
\left(\sum_{j=1}^{n}f'\left(B^{H,K}_{t_{j-1}}\right)1_{(t_{j-1},t_j]}(.)\right)
\underset{n\rightarrow\infty}{\longrightarrow}
f'\left(B^{H,K}_{.}\right)1_{(0,t]}(.) \mbox{ in }
L^2(\Omega;\mathcal{H}).$$ Indeed, using the quasi-helix property of
$B^{H,K}$, we obtain
\begin{eqnarray*} &&E\left|\vert|
\sum _{j=1}^{n}\left[ f'\left(B^{H,K}_{t_{j-1}}\right)-
f'\left(B^{H,K}_{
.}\right)\right] 1_{(t_{j-1},t_j]}(.)\vert| ^{2} _{{\cal{H}}}\right| \\
&=&E \sum_{j,l=1}^{n}\int _{t_{j-1}}^{t_{j}}
\int_{t_{l-1}}^{t_{l}}\left| f'\left(B^{H,K}_{t_{j-1}}\right)-
f'\left(B^{H,K}_{ u}\right)\right|\left|
f'\left(B^{H,K}_{t_{l-1}}\right)- f'\left(B^{H,K}_{
v}\right)\right|\frac{\partial ^{2}R}{\partial u \partial
v}(u,v)dudv\\
&\leq & 2^{1-K}(\sup_{x\in \mathbb{R}}\vert f''(x)\vert)^{2}
\sup_{i=1,\ldots ,n} \left| t_i-t_{i-1}\right|
^{2HK}\sum_{j,l=1}^{n}\int _{t_{j-1}}^{t_{j}}
\int_{t_{l-1}}^{t_{l}}\frac{\partial ^{2}R}{\partial u \partial
v}(u,v)dudv \\
&=& 2^{1-K}(\sup_{x\in \mathbb{R}}\vert f''(x)\vert)^{2}
\sup_{i=1,\ldots ,n} \left| t_i-t_{i-1}\right| ^{2HK}R(t,t) \to
_{n\to \infty }0
\end{eqnarray*}

 This with (\ref{I+J}) and (\ref{J}) implies that
 $I^n_1$ converges   to
$\delta^{B^{H,K}}\left(f'\left(B^{H,K}_{.}\right)1_{(0,t]}(.)\right)$
in $L^2(\Omega)$. Therefore (\ref{ito}) is established due to
(\ref{basic}). \qed \end{proof}

\vskip0.5cm

 The proof of the case $2HK=1$ is based on a preliminary result
 concerning the quadratic variation of the bifractional Brownian
 motion. It was proved in \cite{russo} using the  stochastic
 calculus via regularization.

\begin{lemma}\label{lemma 2hk=1}Suppose that $2HK=1$, then
\begin{eqnarray*}
V_t^n:=\sum_{j=1}^{n}\left(B^{H,K}_{t_{j}}-B^{H,K}_{t_{j-1}}\right)^2\underset{n\rightarrow\infty}{\longrightarrow}
\frac{1}{2^{k-1}}t\ \mbox{ in } L^2(\Omega).
\end{eqnarray*}
\end{lemma}
\begin{proof}
A straightforward calculation shows that,
\begin{eqnarray*}EV_t^{n}=\frac{t}{n}\sum_{j=1}^{n}h(j)+
\frac{t}{2^{k-1}}\underset{n\rightarrow\infty}{\longrightarrow}
\frac{t}{2^{k-1}}.
\end{eqnarray*}
To obtain the conclusion it suffices to show that
\begin{eqnarray*}\lim_{n\rightarrow\infty}E(V_t^{n})^2=(\frac{t}{2^{k-1}})^2.
\end{eqnarray*}
In fact we have,
\begin{eqnarray*}E(V_t^{n})^2=\sum_{i,j=1}^{n}
E\left((B^{H,K}_{t_{i}}-B^{H,K}_{t_{i-1}})(B^{H,K}_{t_{j}}-B^{H,K}_{t_{j-1}})\right)^2
\end{eqnarray*}
Denote by
$$\mu_{n}(i,j)=E\left((B^{H,K}_{t_{i}}-B^{H,K}_{t_{i-1}})(B^{H,K}_{t_{j}}-B^{H,K}_{t_{j-1}})\right)^2$$
It follows by linear regression that
\begin{eqnarray*}\hspace{-3cm}&&\mu_n(i,j)
=E\left(N_1^2\left|\theta_n(i,j)N_1 +\sqrt{\delta_n(i,j)
-\left(\theta_n(i,j)\right)^2}N_2\right|^{2}\right)
\end{eqnarray*}
 where $N_1$ and $N_2$ two independent normal random variables,

\begin{eqnarray*}\theta_n(i,j)&:=&E\left((B^{H,K}_{t_{i}}-B^{H,K}_{t_{i-1}})
(B^{H,K}_{t_{j}}-B^{H,K}_{t_{j-1}})\right)\\&=&\frac{t}{2^Kn}\left[
(i^{2H}+j^{2H})^K-2|j-i|-(i^{2H}+(j-1)^{2H})^K+|j-i-1|\right.\\&&\left.-
((i-1)^{2H}+j^{2H})^K+|j-i+1|+((i-1)^{2H}+(j-1)^{2H})^K\right]\end{eqnarray*}
and
\begin{eqnarray*}\delta_n(i,j)&:=&E\left(B^{H,K}_{t_{i}}-B^{H,K}_{t_{i-1}}\right)^2
E\left(B^{H,K}_{t_{j}}-B^{H,K}_{t_{j-1}}\right)^2.\end{eqnarray*}
Hence
\begin{eqnarray*}\mu_n(i,j)
=2(\theta_n(i,j))^2+\delta_n(i,j)
\end{eqnarray*}
For $1\leq i<j$, we define a function $f_j :
(1,\infty)\rightarrow\R$, by
\begin{eqnarray*}f_j(x)&=&((x-1)^{2H}+j^{2H})^K-((x-1)^{2H}+(j-1)^{2H})^K\\&&-(x^{2H}+j^{2H})^K+(x^{2H}+(j-1)^{2H})^K
\end{eqnarray*} We compute
\begin{eqnarray*}f_j'(x)&=&\left(\frac{(x-1)^{2H}+j^{2H}}{(x-1)^{2H}}\right)^{K-1}
-\left(\frac{(x-1)^{2H}+(j-1)^{2H}}{(x-1)^{2H}}\right)^{K-1}
\\&&
-\left(\frac{x^{2H}+j^{2H}}{x^{2H}}\right)^{K-1}+\left(\frac{x^{2H}+(j-1)^{2H}}{x^{2H}}\right)^{K-1}\\&:=&g(x-1)-g(x)\geq0
\end{eqnarray*}
Hence $f_j$ is increasing and positive, since the  function
$$g(x)=\left(1+\frac{j^{2H}}{x^{2H}}\right)^{K-1}
-\left(1+\frac{(j-1)^{2H}}{x^{2H}}\right)^{K-1}$$ is decreasing on
$(1,\infty)$. This implies that for every $1\leq i<j$
\begin{eqnarray*}|\theta_n(i,j)|=\frac{t}{2^Kn}f_j(i)\leq \frac{t}{2^Kn}f_j(j)\leq\frac{t}{n}|h(j)|\end{eqnarray*}
and $|\theta_n(i,i)|=\frac{t}{n}|h(i)+2|$ for any $i\geq1.$\\Thus
\begin{eqnarray*}\sum_{i,j=1}^{n}\theta_n(i,j)^2\leq\frac{2t^2}{n^2}\sum_{\underset{i,j=1}{i<j}}^{n}h(j)^2+
\frac{t^2}{n^2}\sum_{i=1}^{n}(h(i)+2)^2.
\end{eqnarray*}
Combining this with (\ref{2HK=1}), we obtain that
$\sum_{i,j=1}^{n}\theta_n(i,j)^2$ converges to $0$ as
$n\rightarrow\infty.$\\ On the other hand, by (\ref{2HK=1})
\begin{eqnarray*}\sum_{i,j=1}^{n}\delta_n(i,j)=\frac{t^2}{n^2}\sum_{i,j=1}^{n}\left(h(i)+\frac{1}{2^{K-1}}\right)
\left(h(j)+\frac{1}{2^{K-1}}\right)\underset{n\rightarrow\infty}{\longrightarrow}\left(\frac{t}{2^{K-1}}\right)^2.
\end{eqnarray*}

 Consequently, $E(V_t^n)^2$ converges to
 $\left(\frac{t}{2^{K-1}}\right)^2$ as $n\rightarrow\infty,$ and the
 conclusion follows.
\qed \end{proof}

\vskip1cm

\begin{proof}[Proof of the Theorem \ref{ito unidim} in the case
$2HK=1$.]
 In this case we shall prove that
 \begin{eqnarray}I^n_1
\underset{n\rightarrow\infty}{\longrightarrow}
\int_{0}^{t}f'\left(B^{H,K}_{s}\right)\delta B^{H,K}_{s} \mbox{ in }
L^2(\Omega),
\end{eqnarray}

\begin{eqnarray}I^n_2
\underset{n\rightarrow\infty}{\longrightarrow}
\left(\frac{1}{2}-\frac{1}{2^K}\right)\int_{0}^{t}f''\left(B^{H,K}_s\right)ds
\mbox{ in } L^2(\Omega),\label{I_2}
\end{eqnarray}and
\begin{eqnarray}J^n
\underset{n\rightarrow\infty}{\longrightarrow}
\frac{1}{2^K}\int_{0}^{t}f''\left(B^{H,K}_s\right)ds \mbox{ in }
L^2(\Omega).
\end{eqnarray}
To prove (\ref{I_2}), it is enough to establish that
\begin{eqnarray*}E_n:=\left(E\left|I^n_2-
\left(\frac{1}{2}-\frac{1}{2^K}\right)\sum_{j=1}^{n}f''\left(B^{H,K}_{t_{j-1}}\right)\left(t_{j}-t_{j-1}\right)\right|^2\right)^{1/2}
\underset{n\rightarrow\infty}{\longrightarrow}0.
\end{eqnarray*}
Appliquing the Minkowsky  inequality and (\ref{2HK=1}), we obtain
\begin{eqnarray*}E_n&\leq&C\sum_{j=1}^{n}
 \left|\frac{1}{2^K}(t_j^{2H}+t_{j-1}^{2H})^K
-\frac{1}{2}(t_{j}+t_{j-1})\right|
\\&\leq&\frac{C}{2n}\sum_{j=1}^{n}\left|2j-1-\frac{2}{2^K}(j^{2H}+(j-1)^{2H})^K\right|
\\&=&\frac{C}{2n}\sum_{j=1}^{n}h(j)
\leq\frac{C}{2n}\sum_{j=1}^{n}\frac{1}{j}\underset{n\rightarrow\infty}{\longrightarrow}0.
\end{eqnarray*}
By using Lemma \ref{lemma 2hk=1}, we conclude that
\begin{eqnarray*}
J^n\underset{n\rightarrow\infty}{\longrightarrow}\frac{1}{2^K}\int_0^tf''\left(\overline{B}^{H,K}_{s}\right)ds
\ \mbox{ in } L^2(\Omega).
\end{eqnarray*}
The rest of the proof is same as in the case $2HK>1.$ \qed
\end{proof}

\vskip0.5cm

Let us regard now the Tanaka formula. As in the case of the standard
fractional Brownian motion, it will involve the so-called {\em
weighted local time } $L_t^{^ x}$ ($x\in \mathbb{R}$, $t\in [0,T]) $
 of $B^{H,K}$ defined as the density of the
occupation measure
\begin{eqnarray*}
A\in\mathcal{B}(\R)\longrightarrow2HK\int_0^t1_{A}(B^{H,K}_s)s^{2HK-1}ds.
\end{eqnarray*}

\begin{theorem}Let $\left(B^{H,K}_t,t\in[0,T]\right)$ be a bifractional
Brownian motion with $2HK\geq1$. Then for each $t\in[0,T]$, $x\in\R$
the following formula holds
\begin{eqnarray}\left|B^{H,K}_t-x\right|=|x|+\int_0^tsign(B_s-x)\delta
B^{H,K}_s+L_t^x.
\end{eqnarray}
\end{theorem}
{\bf Proof: }Let $p_{\varepsilon}(y)=\frac{1}{\sqrt{2\pi
\varepsilon}}e^{-\frac{y^2}{2\varepsilon}}$  be the Gaussian kernel
and put
\begin{eqnarray*}F'_{\varepsilon}(z)=2\int_{-\infty}^{z}p_{\varepsilon}(y)dy
- 1,
\end{eqnarray*} and
\begin{eqnarray*}F_{\varepsilon}(z)=\int_{0}^{z}F'_{\varepsilon}(y)dy.
\end{eqnarray*}
By the Theorem \ref{ito unidim} we have
\begin{eqnarray}
F_{\varepsilon}\left(B^{H,K}_{t}-x\right)&=&F_{\varepsilon}\left(-x\right)+\int_{0}^{t}
F'_{\varepsilon}\left(B^{H,K}_{s}-x\right)\delta
B^{H,K}_{s}\nonumber
\\&&+HK
\int_{0}^{t}p_{\varepsilon}\left(B^{H,K}_{s}-x\right)s^{2HK-1}ds.\label{ito
pour F}
\end{eqnarray}
Using (a slightly adaptation of) Proposition 9 in \cite{russo}, one can prove that
\begin{eqnarray}L_t^x=\lim_{\varepsilon\rightarrow0}2HK\int_0^tp_{\varepsilon}(B^{H,K}_s-x)s^{2HK-1}ds\
\mbox{ in } L^2(\Omega).\label{local times }
\end{eqnarray} and $L_t^x$ admits the following chaotic
representation into multiple stochastic integrals (here $I_{n}$
represents the multiple integral with respect to the bifBm)
\begin{eqnarray}\label{Lchaos}L_t^x=2HK\sum_{n=0}^{\infty}\int_0^t\frac{p_{s^{2HK}}(x)}{s^{(n-2)HK+1}}H_n
\left(\frac{x}{s^{HK}}\right)I_n(1_{[0,s]}^{\otimes^n})ds
\end{eqnarray}
where $H_n$ is the nth Hermite polynomial defined as
\begin{eqnarray*}H_n(x)=\frac{(-1)^n}{n!}e^{x^2/2}\frac{d^n}{d
x^n}(e^{-x^2/2})\quad \mbox{ for every } n\geq1.
\end{eqnarray*}
We have  $F_{\varepsilon}(x)\rightarrow |x|$ as
$\varepsilon\rightarrow0$ and since $F_{\varepsilon}(x)\leq |x|$,
then by  Lebesgue's dominated convergence theorem we obtain that
$F_{\varepsilon}(B^{H,K}_t-x)$ converges to $|B^{H,K}_t-x|$ in
$L^2(\Omega)$ as $\varepsilon\rightarrow0$.

On the other hand, since $0\leq F'_{\varepsilon}(x)\leq 1$ and
$F'_{\varepsilon}(x)\rightarrow sign(x)$ as $\varepsilon$ goes to
$0$ the  Lebesgue's dominated convergence theorem in
$L^{2}(\Omega\times[0,T]^2; P\otimes\frac{\partial ^{2}R}{\partial u
\partial v}(u,v)dudv)$ implies that $
F'_{\varepsilon}\left(B^{H,K}_{.}-x\right)$ converges to
$sign\left(B^{H,K}_{.}-x\right)$  in $L^{2}(\Omega; {\cal{H}})$ as
$\varepsilon$ goes to $0$  because
\begin{eqnarray*}
&&E\Vert F'_{\varepsilon}\left(B^{H,K}_{.}-x\right)
-sign\left(B^{H,K}_{.}-x\right) \Vert ^{2} _{|\cal{H}|}\\
&=& E\int_{0}^{T}
\int_{0}^{T}\left|F'_{\varepsilon}\left(B^{H,K}_{u}-x\right)
-sign\left(B^{H,K}_{u}-x\right)\right|\left|F'_{\varepsilon}\left(B^{H,K}_{v}-x\right)
-sign\left(B^{H,K}_{v}-x\right)\right|\\&&\times\frac{\partial
^{2}R}{\partial u
\partial v}(u,v)dudv.
\end{eqnarray*}
Consequently, from the above convergences and (\ref{basic})
\begin{eqnarray*}
\int_{0}^{t} F'_{\varepsilon}\left(B^{H,K}_{s}-x\right)\delta
B^{H,K}_{s}\underset{\varepsilon\rightarrow 0}{\longrightarrow}
\int_{0}^{t}sign\left(B^{H,K}_{s}-x\right)\delta B^{H,K}_{s}\quad
\mbox{ in } L^2(\Omega).
\end{eqnarray*}
Then the conclusion follows. \qed

\vskip0.5cm

\section{Tanaka formula for multidimensional bifractional Brownian
motion}

Given two vectors $H=(H_{1}, \ldots , H_{d})\in [0,1]^{d}$ and
$K=(K_{1}, \ldots , K_{d})\in(0,1] ^{d}$, we introduce the
$d$-dimensional bifractional Brownian motion
$$B^{H,K} =\left(B^{H_1,K_1},...,B^{H_d,K_d}\right)$$
as a centered Gaussian vector whose component are independent
one-dimensional bifractional Brownian motions.

\vskip0.3cm

We extend the It\^o formula to the multidensional case.

\begin{theorem}Let $B^{H,K}=\left(B^{H_1,K_1},...,B^{H_d,K_d}\right)$
be a $d$-dimensional bifractional Brownian motion, and let $f$ be a
function of class
$\mathrm{C}^2\left(\mathbb{R}^d,\mathbb{R}\right)$. We assume that
$2H_iK_i>1$ for any $i=1,...n$, then
\begin{eqnarray}
f\left(B^{H,K}_{t}\right)=f\left(0\right)+\sum_{i=1}^{d}\int_{0}^{t}\frac{\partial
f}{\partial x_i}\left(B^{H_i,K_i}_{s}\right)\delta
B^{H_i,K_i}_{s}+\sum_{i=1}^{d} H_iK_i
\int_{0}^{t}\frac{\partial^2f}{\partial
x_i^2}\left(B^{H,K}_{s}\right)s^{2H_iK_i-1}ds.\label{ito2}
\end{eqnarray}
\end{theorem}
\begin{proof}Let us fix $t>0$ and a partition
$\{t_j=\frac{jt}{n}; j=0,...,n\}$ of $[0,t]$.  As in above we may
assume  that f has compact support, and so $f$, $f'$ and $f''$ are
bounded. Using Taylor expansion, we have
\begin{eqnarray*}f\left(B^{H,K}_t\right)&=&f(0)+\sum_{j=1}^{n}\sum_{i=1}^{d}\frac{\partial
f}{\partial x_i}\left(B^{H,K}_{t_{j-1}}\right)
\left(B^{H_i,K_i}_{t_{j}}-B^{H_i,K_i}_{t_{j-1}}\right)\\&&+\frac{1}{2}\sum_{j=1}^{n}
\sum_{i,l=1}^{d}\frac{\partial^2f}{\partial x_i\partial
x_l}\left(\overline{B}^{H,K}_{j}\right)
\left(B^{H_i,K_i}_{t_{j}}-B^{H_i,K_i}_{t_{j-1}}\right)\left(B^{H_l,K_l}_{t_{j}}-B^{H_l,K_l}_{t_{j-1}}\right)
\\&:=&f(0)+
I^n + J^n.
\end{eqnarray*}
where
$\overline{B}^{H,K}_{j}=B^{H,K}_{t_{j-1}}+\theta_j\left(B^{H,K}_{t_{j}}-B^{H,K}_{t_{j-1}}\right)$,
and $\theta_j$ is a random variable in $(0,1).$\\
We show that $J^n$ converges to 0 in $L^2(\Omega)$ as
$n\rightarrow\infty$.
\begin{eqnarray*}E|J^n|^2&\leq&
n/4\sum_{j=1}^{n}E\left(\sum_{i,l=1}^{d}\frac{\partial^2f}{\partial
x_i\partial x_l}\left(\overline{B}^{H,K}_{j}\right)
\left(B^{H_i,K_i}_{t_{j}}-B^{H_i,K_i}_{t_{j-1}}\right)\left(B^{H_l,K_l}_{t_{j}}-B^{H_l,K_l}_{t_{j-1}}\right)\right)^2\\
&\leq&\frac{Cd^2}{4}n\sum_{j=1}^{n}\sum_{i,l=1}^{d}
E\left(B^{H_i,K_i}_{t_{j}}-B^{H_i,K_i}_{t_{j-1}}\right)^2E\left(B^{H_l,K_l}_{t_{j}}-B^{H_l,K_l}_{t_{j-1}}\right)^2\\
&\leq&Cd^2n\sum_{i,l=1}^{d}\sum_{j=1}^{n}|t_j-t_{j-1}|^{2(H_iK_i+H_lK_l)}=Cd^2\sum_{i,l=1}^{d}\frac{t^{4HK}}
{n^{2(H_iK_i+H_lK_l-1)}}
\underset{n\rightarrow\infty}{\longrightarrow}0
\end{eqnarray*}
According to (\ref{ip}), we get
\begin{eqnarray*}I^n&=&\sum_{j=1}^{n}\sum_{i=1}^{d}\frac{\partial
f}{\partial
x_i}\left(B^{H,K}_{t_{j-1}}\right)\left(\delta^{B^{H_i,K_i}}(1_{(t_{j-1},t_j]})\right)\\
&=&\sum_{i=1}^{d}\left[\delta^{B^{H_i,K_i}}\left(\sum_{j=1}^{n}\frac{\partial
f}{\partial
x_i}\left(B^{H,K}_{t_{j-1}}\right)1_{(t_{j-1},t_j]}(.)\right)\right.\\&&+\left.
\sum_{j=1}^{n}\frac{\partial^2 f}{\partial
x_i^2}\left(B^{H,K}_{t_{j-1}}\right)\langle1_{(0,t_{j-1}]},1_{(t_{j-1},t_j]}\rangle_{\mathcal{H}}\right]\\&=&
\sum_{i=1}^{d}\left[I_1^{n,i}+I_2^{n,i} \right]\end{eqnarray*} As
the similar way in the above theorem, we obtain that for every
$i=1,...,d$
$$I_2^{n,i}\underset{n\rightarrow\infty}{\longrightarrow}
 H_iK_i \int_{0}^{t}\frac{\partial^2f}{\partial
x_i^2}\left(B^{H,K}_{s}\right)s^{2H_iK_i-1}ds\mbox { in
}L^2(\Omega).$$ We show that for every $i\in\{1,...,d\}$
$$I_1^{n,i}\underset{n\rightarrow\infty}{\longrightarrow}
 \int_{0}^{t}\frac{\partial f}{\partial
x_i}\left(B^{H,K}_{s}\right)\delta B^{H_i,K_i}_{s}\mbox { in
}L^2(\Omega).$$ We set
\begin{eqnarray*}u_s^{n,i}=\sum_{j=1}^{n}\frac{\partial
f}{\partial
x_i}\left(B^{H,K}_{t_{j-1}}\right)1_{(t_{j-1},t_j]}(s)-\frac{\partial
f}{\partial x_i}\left(B^{H,K}_{s}\right)1_{(0,t]}(s).
\end{eqnarray*}
By inequality (\ref{ineg1}), we have
\begin{eqnarray*}E\left(\delta^{B^{H_i,K_i}}(u^{n,i})\right)^2
\leq E\|u^{n,i}\|_{\left| {\cal{H}}^i\right|}^2+E\Vert Du^{n,i}\Vert
^{2} _{ \left| {\cal{H}}^i\right|\otimes \left| {\cal{H}}^i\right| }
\end{eqnarray*} where ${\cal{H}}^i$ is the Hilbert space associated
to $B^{H_i,K_i}$ and $R_i$ its covariance function.\\For every $r,
s\leq t$
\begin{eqnarray*}D_ru_s^{n,i}=\sum_{j=1}^{n}\frac{\partial^2
f}{\partial
x_i^2}\left(B^{H,K}_{t_{j-1}}\right)1_{(0,t_{j-1}]}(r)1_{(t_{j-1},t_j]}(s)-\frac{\partial^2
f}{\partial x_i^2}(B_s^{H,K}) 1_{(0,s]}(r)\end{eqnarray*} we remark
that $D_ru_s^{n,i}$ and $u_s^{n,i}$ converge to zero as
$n\longrightarrow\infty$ for any $r, s\leq t$. Since the first and
second partial derivatives of $f$ are bounded, then by using the
Lebesgue dominated convergence theorem and the expression of the
norm $\left| {\cal{H}}^i\right|\otimes \left| {\cal{H}}^i\right| $
we obtain that
\begin{eqnarray*}\delta^{B^{H_i,K_i}}(u^{n,i})\underset{n\rightarrow\infty}{\longrightarrow}0
 \mbox { in
}L^2(\Omega).
\end{eqnarray*}
The proof is thus complete.

\qed \end{proof}

\vskip0.5cm

One can easily generalize the above theorem to the case when the
function $f$ depends on time.
\begin{theorem}
\label{timedep} Let $f\in \mathrm{C}^{1,2}\left([0,T]
\times\mathbb{R}^d,\mathbb{R}\right)$ and
$B^{H,K}=\left(B^{H_1,K_1},...,B^{H_d,K_d}\right)$ be a
$d$-dimensional bifBm with  $2H_iK_i>1$ for any $i=1,...,n$. Then
\begin{eqnarray}
f\left(t, B^{H,K}_{t}\right)&=& f\left(0, 0\right)+\int_{0}^{t}
\frac{\partial f}{\partial s}(s,B^{H,K}_{s})ds
 + \sum_{i=1}^{d}\int_{0}^{t}\frac{\partial f}{\partial
x_i}\left(s, B^{H,K}_{s}\right)\delta
B^{H_i,K_i}_{s} \nonumber \\\
&&+\sum_{i=1}^{d} H_iK_i \int_{0}^{t}\frac{\partial^2f}{\partial
x_i^2}\left(s, B^{H,K}_{s}\right)s^{2H_iK_i-1}ds.\label{ito3}
\end{eqnarray}

\end{theorem}

\vskip0.5cm

We consider twice of the kernel of the $d-$dimensional Newtonian
potential

$$U(z)=\left\{
\begin{array}{lr}
-\frac{\Gamma(d/2-1)}{2\pi^{d/2}}\frac{1}{|z|^{d-2}}
\quad\mbox{if}\ d\geq3\\
\frac{1}{\pi}log|z|\qquad\qquad\mbox{if}\ d=2. \end{array} \right.$$
 Set
\begin{equation}
\label{ubar} \bar{U}(s,z)=\frac{1}{\prod_{j=1}^{d}\sqrt{2H_jK_j}}
s^{\theta}U\left(\frac{(z_1-x_1)}{\sqrt{2H_1K_1}}s^{1/2-H_1K_1},...,\frac{(z_d-x_d)}{\sqrt{2H_dK_d}}s^{1/2-H_dK_d}\right)
\end{equation}
where $x=(x_1,...,x_d)\in\R^d$ and
$0<\gamma:=\frac{1}{2}(2-d)+\theta+(d-2)(HK)^*-\sum_{i=1}^{d}H_iK_i$
 with $(HK)^*~=~max\{H_1K_1,\ldots ,H_dK_d\}$.
\vskip0.5cm

We shall prove the following Tanaka formula. It will involve a
multidimensional weighted local time  which is an extension of the
one-dimensional local time  given by (\ref{Lchaos}). Note for any
dimension $d\geq 2$ the local time is not a random variable anymore
and it is a distribution in the Watanabe's sense.
\begin{theorem}
\label{multi-tanaka} Let $\bar{U}$  as above and let $B^{H,K}=\left(B^{H_1,K_1},...,B^{H_d,K_d}\right)$ be a $d$-dimensional
bifBm with  $2H_iK_i>1$ for any $i=1,...d$. Then the following formula holds in the Watanabe space $\mathbb{D}_2^{\alpha-1 }$
for any $\alpha<\frac{1}{2(HK)^*}-d/2$.

\begin{eqnarray}\bar{U}(t,B_t^{H,K})=\bar{U}(0,0)+\int_{0}^{t}\partial_s\bar{U}(s,B_s^{H,K})ds+
\sum_{i=1}^{d} \int_{0}^{t}\frac{\partial \bar{U}(s,B_s^{H,K})}{\partial x_{i}} \delta
B_s^{H_i,K_i}+L^{\theta}(t,x)\label{tanaka-multidim}
\end{eqnarray}
where the generalized weighted local time $L^{\theta}(t,x)$ is
defined as
$$L^{\theta}(t,x)=\sum_{n=(n_1,...,n_d)
}\int_0^t\prod_{i=1}^{d}\frac{p_{s^{2H_iK_i}}(x_i)}
{s^{\frac{1}{2}+(n_i-1)H_iK_i}}
H_{n_i}\left(\frac{x_i}{\sqrt{s^{2H_iK_i}}}\right)I_{n_i}^i(1_{[0,s]}^{\otimes^{n_i}})s^{\theta}ds.$$\\
\end{theorem}
{\bf Proof: } We regularize the function $\bar{U}$ by standard
convolution. Put $\bar{U}_{\varepsilon}=p_{\varepsilon}^d*\bar{U}$,
with $p_{\varepsilon}^d$ is the Gaussian kernel on $\R^d$ given by
$$p_{\varepsilon}^d(x)=\prod_{i=1}^{d}p_{\varepsilon}(x_i)=\prod_{i=1}^{d}\frac{1}{\sqrt{2\pi
\varepsilon}}e^{-\frac{x^2}{2\varepsilon}}, \quad\forall
x=(x_1,...,x_d)\in\R^d.$$ Using the above It\^o formula we have
\begin{eqnarray*}
\bar{U}_{\varepsilon}\left(t,B^{H,K}_{t}\right)&=&\bar{U}_{\varepsilon}\left(0,0\right)+\int_{0}^{t}\frac{\partial
\bar{U}_{\varepsilon}}{\partial
s}\left(s,B^{H,K}_{s}\right)ds+\sum_{i=1}^{d}\int_{0}^{t}\frac{\partial
\bar{U}_{\varepsilon}}{\partial x_i}\left(s,B^{H,K}_{s}\right)\delta
B^{H_i,K_i}_{s}\\&&+\sum_{i=1}^{d} H_iK_i
\int_{0}^{t}\frac{\partial^2 \bar{U}_{\varepsilon}}{\partial
x_i^2}\left(s,B^{H,K}_{s}\right)s^{2H_iK_i-1}ds.\\&=&\bar{U}_{\varepsilon}\left(0,0\right)+I_1^{\varepsilon}(t)+I_2^{\varepsilon}(t).
\end{eqnarray*}
On the other hand, if $V(z)=U(a_1z_1,...,a_dz_d)$ and
$V_{\varepsilon}=p^d_{\varepsilon}*V$ we have
$$\frac{1}{2}\sum_{i=1}^{d}\frac{1}{a_i^2}\frac{\partial^2V_{\varepsilon}}{\partial
z_i^2}(z)=p^d_{\varepsilon}(a_1z_1,...,a_dz_d).$$
Hence\begin{eqnarray*}
I_2^{\varepsilon}(t)=\frac{1}{\prod_{j=1}^{d}\sqrt{2H_jK_j}}\int_0^t
p^d_{\varepsilon}\left(c_1(s)(B_s^{H_1,K_1}-x_1),...,c_d(s)(B_s^{H_d,K_d}-x_d)\right)s^{\theta}ds
\end{eqnarray*}
where $c_i(s)=\frac{s^{1/2-H_iK_i}}{\sqrt{2H_iK_i}}$ for every $i=
1,...,d$. The next step is to find the chaotic expansion of the last
term $I_{2}^{\varepsilon}$. By Strook formula, we have
\begin{eqnarray*}p_{\varepsilon}\left(c_i(s)(B_s^{H_i,K_i}-x_i)\right)=\sum_{n=0}^{\infty}\frac{1}{n!}
I_n^i\left(ED^n_.p_{\varepsilon}\left(c_i(s)(B_s^{H_i,K_i}-x_i)\right)\right).
\end{eqnarray*}
and
\begin{eqnarray*}&&ED^n_.p_{\varepsilon}\left(c_i(s)(B_s^{H_i,K_i}-x_i)\right)=c_i(s)^nEp_{\varepsilon}^{(n)}\left(c_i(s)(B_s^{H_i,K_i}-x_i)\right)
1_{[0,s]}^{\otimes^n}(.)\\&=&c_i(s)^nn!(\frac{\varepsilon}{c_i(s)^2}+s^{2H_iK_i})^{-n/2}p_{s^{2H_iK_i}+\frac{\varepsilon}{c_i(s)^2}}(x_i)
\frac{H_n\left(\frac{x_i}{\sqrt{s^{2H_iK_i}+\frac{\varepsilon}{c_i(s)^2}}}\right)}{c_i(s)^{n+1}}1_{[0,s]}^{\otimes^n}(.)
\\&=&\frac{n!}{c_i(s)}(\frac{\varepsilon}{c_i(s)^2}+s^{2H_iK_i})^{-n/2}p_{s^{2H_iK_i}+\frac{\varepsilon}{c_i(s)^2}}(x_i)
H_n\left(\frac{x_i}{\sqrt{s^{2H_iK_i}+\frac{\varepsilon}{c_i(s)^2}}}\right)1_{[0,s]}^{\otimes^n}(.)\\&:=&\frac{n!}{c_i(s)}
\beta_{n,\varepsilon}^{i}(s)1_{[0,s]}^{\otimes^n}(.)\end{eqnarray*}
Consequently
\begin{eqnarray*}p^d_{\varepsilon}\left(c_1(s)(B_s^{H_1,K_1}-x_1),...,c_d(s)(B_s^{H_d,K_d}-x_d)\right)=
\sum_{n=(n_1,...,n_d)\in\mathbb{N}^d}\prod_{i=1}^{d}\frac{\beta_{n_i,\varepsilon}^{i}(s)}{c_i(s)}I_{n_i}^i(1_{[0,s]}^{\otimes^{n_i}})
\end{eqnarray*}
and that
\begin{eqnarray*}I_2^{\varepsilon}(t)
&=&\sum_{n=(n_1,...,n_d)
\in\mathbb{N}^d}\int_0^t\prod_{i=1}^{d}\frac{\beta_{n_i,\varepsilon}^{i}(s)}
{s^{\frac{1}{2}-H_iK_i}}I_{n_i}^i(1_{[0,s]}^{\otimes^{n_i}})s^{\theta}ds\\
&=&
\sum_{n=(n_1,...,n_d)}\int_0^t\prod_{i=1}^{d}\frac{p_{s^{2H_iK_i}+\frac{\varepsilon}{c_i(s)^2}}(x_i)}
{(\frac{\varepsilon}{c_i(s)^2}+s^{2H_iK_i})^{n_i/2}s^{\frac{1}{2}-H_iK_i}}
H_{n_i}\left(\frac{x_i}{\sqrt{s^{2H_iK_i}+\frac{\varepsilon}{c_i(s)^2}}}\right)I_{n_i}^i(1_{[0,s]}^{\otimes^{n_i}})s^{\theta}ds
\end{eqnarray*}
This term (in fact, slightly modified) appeared in some other papers
such as Proposition 12 in \cite{ELSTV}, or in \cite{Xiao}. Using
standard arguments we obtain that the last term converges  in
$\mathbb{D}^{\alpha}_2$ to $L^{\theta}(t,x)$ as $\varepsilon$ goes
to $0$, with $\alpha<\frac{1}{2(HK)^{*}}-d/2$.

\vskip0.5cm

The rest of the proof is to show that the following convergences are
holds:\\
 For every $i=1,...,d$
\begin{eqnarray}\int_{0}^{t}\partial_i \bar{U}_{\varepsilon}\left(s,B^{H,K}_{s}\right)\delta
B^{H_i,K_i}_{s}\overset{\mathbb{D}^{\alpha-1}_2}{\underset{\varepsilon\rightarrow0}{\longrightarrow}}\int_{0}^{t}\partial_i
\bar{U}\left(s,B^{H,K}_{s}\right)\delta
B^{H_i,K_i}_{s}.\label{converg. of garadient in D}
\end{eqnarray}
\begin{eqnarray}
\int_0^t\partial_s\bar{U}_{\varepsilon}(s,B_s^{H,K})ds\overset{\mathbb{D}^{\alpha}_2}{\underset{\varepsilon\rightarrow0}{\longrightarrow}}
\int_0^t\partial_s\bar{U}(s,B_s^{H,K})ds
\label{converg.ofpartial_sUin D}
\end{eqnarray} and
\begin{eqnarray}\bar{U}_{\varepsilon}(t,B_t^{H,K})\overset{\mathbb{D}^{\alpha}_2}{\underset{\varepsilon\rightarrow0}{\longrightarrow}}
\bar{U}(t,B_t^{H,K}). \label{converg. of U in D}
\end{eqnarray}
We start with the convergence (\ref{converg. of garadient in D}).
Fix $i\in\{1,...,d\}$, we note $g_{\varepsilon}^i(s,z)=\partial_i
\bar{U}_{\varepsilon}\left(s,z\right)$. By the formal relation
($\delta$ is the Dirac distribution)
$$\int_{\mathbb{R^{d}}} f(y)\delta (x-y)dy =f(x)$$
 we can write (this is true in the sense of Watanabe
distributions)
$$g_{\varepsilon}^i(s,B^{H,K}_{s})=\int_{\R^d}g_{\varepsilon}^i(s,y)\delta\left(B^{H,K}_{s}-y\right)dy.$$
Furthermore (see \cite{ELSTV} but it can be also derived from a
general formula in \cite{kuo})
\begin{eqnarray*}\delta\left(B^{H,K}_{s}-y\right)&=&\prod_{i=1}^{d}\delta\left(B^{H_i,K_i}_{s}-y_i\right)\\&=&
\prod_{i=1}^{d}\left(\sum_{n\geq0}\frac{1}{(R_i(s))^{n/2}}p_{R_i(s)}(y_i)H_n(\frac{y_i}{R_i(s)^{1/2}})I^{i}_{n}\left(
1_{[0,s]}^{\otimes^n}\right)\right)\\&=&
\sum_{n=(n_1,...,n_d)}A_{n}(s,y)I_n(1_{[0,s]}^{\otimes^{|n|}})
\end{eqnarray*}
where $R_{i}(s)=R_{i}(s,s)=s^{2H_i K_i}$,
$A_{n}(s,y)=\prod_{i=1}^{d}\frac{1}{(R_i(s))^{n_i/2}}p_{R_i(s)}(y_i)
H_{n_i}(\frac{y_i}{R_i(s)^{1/2}})$ and
$I_n(1_{[0,s]}^{\otimes^{|n|}}):=\prod_{i=1}^{d}I^{i}_{n_i}\left(
1_{[0,s]}^{\otimes^{n_i}}\right)$ for every $n=(n_1,...,n_d)$.\\
Hence
\begin{eqnarray*}g_{\varepsilon}^i(s,B^{H,K}_s)&=&\sum_{n=(n_1,...,n_d)}\int_{\R^d}g_{\varepsilon}^i(s,y)A_{n}(s,y)dyI_n(1_{[0,s]}^{\otimes^{|n|}})
\\&:=&\sum_{n=(n_1,...,n_d)}B_n^{\varepsilon,i}(s)I_n(1_{[0,s]}^{\otimes^{|n|}})\end{eqnarray*}
and using the chaotic form of the divergence integral
 \begin{eqnarray*}
J_i^{\varepsilon}(t)&:=&\int_0^tg_{\varepsilon}^i\left(s,B^{H,K}_{s}\right)
\delta
B^{H_i,K_i}_{s}\\&=&\sum_{n=(n_1,...,n_d)}I_{n_i+1}^{i}\left[B_n^{\varepsilon,i}(s)1_{[0,s]}^{\otimes^{n_i}}(s_1,...,s_{n_i
})1_{[0,t]}(s)\prod_{\overset{j=1}{j\neq
i}}^{d}I^j_{n_j}\left(1_{[0,s]}^{\otimes^{n_j}}\right)\right] ^{(s)}
\\&=&\sum_{n=(n_1,...,n_d)}I_{n_i+1}^{i}\left[f_{i,n}^{\varepsilon,t}(s_1,...,s_{n_i
},s)\right]\\&=&\sum_{n_i\geq0}I_{n_i+1}^{i}\left[\sum_{n=(n_1,...,\hat{n}_i,...,n_d)}f_{i,n}^{\varepsilon,t}(s_1,...,s_{n_i
},s)\right]
\end{eqnarray*}
where the superscript $(s)$ denoted the symmetrization with respect
to $s_{1}, \ldots , s_{n_{i}}, s$, and
$$f_{i,n}^{\varepsilon,t}(s_1,...,s_{n_i+1
})=\sum_{l=1}^{n_i+1}\frac{1}{n_i+1}B_n^{\varepsilon,i}(s_l)1_{[0,s_l]}^{\otimes^{n_i}}(s_1,...,\widehat{s_l},...,s_{n_i+1
})1_{[0,t]}(s_l)\prod_{\overset{j=1}{j\neq
i}}^{d}I^j_{n_j}\left(1_{[0,s_l]}^{\otimes^{n_j}}\right).$$

Observe here that, since  the components of the vector $B^{H,K}$ are independent, the term $\prod_{\overset{j=1}{j\neq
i}}^{d}I^j_{n_j}\left(1_{[0,s_l]}^{\otimes^{n_j}}\right)$ is viewed as a deterministic function for the integral
$I^{i}_{n_{i}}$. The convergence (\ref{converg. of garadient in D}) is satisfied if the conditions i) and ii) of Lemma 3 in
\cite{ELSTV}  hold. It is easy to verify the condition i), we will prove only the condition
ii).\\
Fixing $i\in\{1,...,d\}$, we can write,
\begin{eqnarray*}\|J_i^{\varepsilon}(t)\|_{\alpha-1,2}^2
&\leq&\sum_{m\geq1}(m+1)^{\alpha-1}\sum_{|n|=n_1+...+n_d=m-1}(n_i+1)!E\left\|f_{i,n}^{\varepsilon,t}\right
\|^2_{{\cal{H}}^{\otimes^{n_i+1}}}\\
&=&\sum_{m\geq1}(m+1)^{\alpha-1}\sum_{|n|=n_1+...+n_d=m-1}(n_i+1)!\\&\times&\int_{[0,T]^{n_i+1}}\int_{[0,T]^{n_i+1}}
\sum_{l,k=1}^{n_i+1}\frac{1}{(n_i+1)^2}|B_n^{\varepsilon,i}(s_l)||B_n^{\varepsilon,i}(r_k)|1_{[0,t]}(s_l)1_{[0,t]}(r_k)
\\&\times&1_{[0,s_l]}^{\otimes^{n_i}}
(s_1,...,\widehat{s_l},...,s_{n_i+1
})1_{[0,r_k]}^{\otimes^{n_i}}(r_1,...,\widehat{r_l},...,r_{n_i+1})\prod_{\underset{j=1}{j\neq
i}}^{d}n_j!R_j(s_l,r_k)^{n_j}\\&\times&
\prod_{q=1}^{n_i+1}\frac{\partial^2R_i}{d s_qd
r_q}(s_q,r_q)dr_1....dr_{n_i+1}ds_1....ds_{n_i+1}\end{eqnarray*}
Since $\alpha<0$ then $(m+2)^{\alpha-1}\leq(m+1)^{\alpha-1}$ and
$(n_i+1)\leq(m+1)$. This implies that
\begin{eqnarray*}\|J_i^{\varepsilon}(t)\|_{\alpha-1,2}^2&\leq&
\sum_{m\geq0}(m+1)^{\alpha}\sum_{|n|=n_1+...+n_d=m}(n_i)!\left[(1-\frac{1}{(n_i+1)})
\right.\\&\times&\left.\int_{[0,T]^{2}}\int_{[0,T]^{2}}|B_n^{\varepsilon,i}(s_1)||B_n^{\varepsilon,i}(r_2)|1_{[0,t]}(s_1)1_{[0,t]}(r_2)
R_i(s_1,r_2)^{n_i-1}\right.\\&\times&\left.\prod_{\underset{j=1}{j\neq
i}}^{d}n_j!R_j(s_1,r_2)^{n_j}1_{[0,s_1]}(s_2)1_{[0,r_2]}(r_1)
\frac{\partial^2R_i}{d s_1d r_1}(s_1,r_1)\frac{\partial^2R_i}{d s_2d
r_2}(s_2,r_2)ds_1ds_2dr_1dr_2\right.\\&+&\left.\frac{1}{(n_i+1)}
\int_{[0,T]^{2}}\int_{[0,T]^{2}}|B_n^{\varepsilon,i}(s_1)||B_n^{\varepsilon,i}(r_1)|1_{[0,t]}(s_1)1_{[0,t]}(r_1)
R_i(s_1,r_1)^{n_i-1}\right.\\&\times&\left.\prod_{\underset{j=1}{j\neq
i}}^{d}n_j!R_j(s_1,r_1)^{n_j}1_{[0,s_1]}(s_2)1_{[0,r_1]}(r_2)
\frac{\partial^2R_i}{d s_1d r_1}(s_1,r_1)\frac{\partial^2R_i}{d s_2d
r_2}(s_2,r_2)ds_1ds_2dr_1dr_2\right]\end{eqnarray*} By integration
we obtain
\begin{eqnarray*}\|J_i^{\varepsilon}(t)\|_{\alpha-1,2}^2&\leq&
\sum_{m\geq0}(m+1)^{\alpha}\sum_{|n|=n_1+...+n_d=m}(n_i)!\left[(1-\frac{1}{(n_i+1)})\right.\\&\times&\left.
\int_{[0,t]^{2}}|B_n^{\varepsilon,i}(s_1)||B_n^{\varepsilon,i}(r_2)|
R_i(s_1,r_2)^{n_i-1}\right.\\&\times&\left.\prod_{\underset{j=1}{j\neq
i}}^{d}n_j!R_j(s_1,r_2)^{n_j} \frac{\partial R_i}{d
s_1}(s_1,r_2)\frac{\partial R_i}{d
r_2}(s_1,r_2)ds_1dr_2\right.\\&+&\left.\frac{1}{(n_i+1)}
\int_{[0,t]^{2}}|B_n^{\varepsilon,i}(s_1)||B_n^{\varepsilon,i}(r_1)|
R_i(s_1,r_1)^{n_i}\right.\\&\times&\left.\prod_{\underset{j=1}{j\neq
i}}^{d}n_j!R_j(s_1,r_1)^{n_j} \frac{\partial^2R_i}{d s_1d
r_1}(s_1,r_1)ds_1dr_1\right].
\end{eqnarray*}
We have for any $1/4\leq\beta\leq1/2$
\begin{eqnarray*}B_n^{\varepsilon,i}(s)&=&\int_{\R^d}g_{\varepsilon}^i(s,y)
A_n(s,y)dy\\&=&\int_{\R^d}g_{\varepsilon}^i(s,y)\prod_{j=1}^{d}
H_{n_j}\left(\frac{y_j}{\sqrt{R_j(s)}}\right)e^{\frac{-\beta
y_j^2}{R_j(s)}}\frac{1}{\sqrt{R_j(s)}^{n_j}}\frac{e^{-(\frac{1}{2}-\beta)\frac{
y_j^2}{R_j(s)}}}{\sqrt{2\pi R_j(s)}}
dy\\&=&\int_{\R^d}p_{\varepsilon}^d(z)dz\int_{\R^d}\partial_i\bar{U}\left(s,y-z\right)\prod_{j=1}^{d}
H_{n_j}\left(\frac{y_j}{\sqrt{R_j(s)}}\right)\frac{e^{\frac{-\beta
y_j^2}{R_j(s)}}}{\sqrt{2\pi }}\frac{e^{-(\frac{1}{2}-\beta)\frac{
y_j^2}{R_j(s)}}}{\sqrt{R_j(s)}^{n_j+1}} dy.
\end{eqnarray*}
Since (see Lemma 11 in \cite{ELSTV})
$$\sup_{z\in\R^d}\prod_{j=1}^{d}|H_{n_j}\left(z_j\right)|e^{-\beta
z_j^2}\leq C
\prod_{j=1}^{d}\frac{1}{\sqrt{n_j!}(n_j\vee1)^{\frac{8\beta-1}{12}}}.$$
and for every $(s,z)\in(0,T]\times\R^d$
\begin{eqnarray*}\left|\partial_i
\bar{U}(s,z)\right|\leq C
s^{\frac{1}{2}(1-d)+\theta}\left|\left((z_1-x_1)s^{-H_1K1},...,(z_d-x_d)s^{-H_dKd}\right)\right|^{1-d}.
 \end{eqnarray*}
Then, for any $s\in(0,T]$
\begin{eqnarray*}|B_n^{\varepsilon,i}(s)|&\leq&C
\prod_{j=1}^{d}\frac{1}{\sqrt{n_j!}(n_j\vee1)^{\frac{8\beta-1}{12}}}
\frac{1}{s^{n_jH_jK_j}}\int_{\R^d}p_{\varepsilon}^d(z)dz\\&\times&
\int_{\R^d}s^{\frac{1}{2}(1-d)+\theta-\sum_{j=1}^{d}H_jK_j}
\frac{e^{-(\frac{1}{2}-\beta)|s^{-HK}y|^2}}{|s^{-HK}(y-z-x)|^{d-1}}dy\\
&\leq&C s^{\gamma-\frac{1}{2}}
\prod_{j=1}^{d}\frac{1}{\sqrt{n_j!}(n_j\vee1)^{\frac{8\beta-1}{12}}}
\frac{1}{s^{n_jH_jK_j}}\int_{\R^d}p_{\varepsilon}^d(z)dz\\&\times&
\int_{\R^d}
\frac{e^{-(\frac{1}{2}-\beta)|s^{-HK}y|^2}}{|(y-z-x)|^{d-1}}dy
\end{eqnarray*}where $s^{- HK}y:=(s^{- H_1K_1}y_1,...,s^{-
H_dK_d}y_d)$§.\\
Let $\eta$ a positive constant such that, for every $s\in(0,T]$,
$j\in\{1,\ldots,d\}$ we have $s^{-2H_jK_j}>\eta.$\\ Combining this
with for any $a,b\in\R$, $a^2\geq\frac{1}{2}(a-b)^2-b^2$, we obtain
\begin{eqnarray*}\int_{\R^d}p_{\varepsilon}^d(z)dz\int_{\R^d}
\frac{e^{-(\frac{1}{2}-\beta)|s^{-HK}y|^2}}{|(y-z-x)|^{d-1}}dy&\leq&
\int_{\R^d}p_{\varepsilon}^d(z)dz\int_{\R^d}
\frac{e^{-\eta(\frac{1}{2}-\beta)|y|^2}}{|(y-z-x)|^{d-1}}dy\\&\leq&
\int_{\R^d}p_{\varepsilon}^d(z)e^{\eta(\frac{1}{2}-\beta)|z+x|^2}dz\int_{\R^d}
\frac{e^{-\frac{\eta}{2}(\frac{1}{2}-\beta)|y-(z+x)|^2}}{|(y-(z+x))|^{d-1}}dy\\&\leq&
Ce^{2\eta(\frac{1}{2}-\beta)|x|^2}\int_{\R^d}p_{\varepsilon}^d(z)e^{2\eta(\frac{1}{2}-\beta)|z|^2}dz
\\&\leq&Ce^{2\eta(\frac{1}{2}-\beta)|x|^2}\int_{\R^d}
\frac{e^{-\frac{|v|^2}{2}}}{\sqrt{2\pi}^d}e^{2\eta(\frac{1}{2}-\beta)\varepsilon|v|^2}dv
\\&\leq&Ce^{2\eta(\frac{1}{2}-\beta)|x|^2}\int_{\R^d}
\frac{e^{-\beta|v|^2}}{{\sqrt{2\pi}^d}}dv<\infty
\end{eqnarray*} since $2\eta\varepsilon\leq 1$ when $\varepsilon $
close to 0.\\ Thus
\begin{eqnarray*}|B_n^{\varepsilon,i}(s)|
\leq C s^{\gamma-\frac{1}{2}}
\prod_{j=1}^{d}\frac{1}{\sqrt{n_j!}(n_j\vee1)^{\frac{8\beta-1}{12}}}
\frac{1}{s^{n_jH_jK_j}}
\end{eqnarray*} where $C$ is a constant depending only on $d,\ H,\ K,\ T,\
x$ and $\beta$.
\\Using this inequality $a^2+b^2\leq2ab$, for every $a, b\in\R_+$, we
conclude  that there exist a constant $C(H,K)$ positive, such that
for every $i=1,\ldots,d$ $$\left|\frac{\partial R_i}{d
r}(r,s)\frac{\partial R_i}{d s}(r,s)\right|\leq C(H,K)
(rs)^{2H_iK_i-1},$$
$$\left|\frac{\partial^2R_i}{d rd s}(r,s)\right|\leq C(H,K) (rs)^{H_iK_i-1}$$and
$$\left|\frac{R_i(r,s)}{(rs)^{H_iK_i}}\right|\leq C(H,K).$$
It follows by anterior inequalities that
\begin{eqnarray*}
\|J_i^{\varepsilon}(t)\|_{\alpha-1,2}^2&\leq&C\sum_{m\geq0}(1+m)^{\alpha}\sum_{|n|=n_1+...+n_d=m}\prod_{j=1}^{d}\frac{1}{(n_j\vee1)^{\frac{8\beta-1}{6}}}
\\&\times&
\int_{[0,t]^{2}}\frac{R_i(r,s)^{n_i-1}}{(rs)^{(n_i-1)H_iK_i}}(rs)^{\gamma-\frac{3}{2}+H_iK_i}\prod_{\underset{j\neq
i}{j=1}}^{d}\frac{R_j(r,s)^{n_j}}{(rs)^{n_jH_jK_j}} drds.
\end{eqnarray*}
We use the selfsimilarity of the covariance kernel $R(r,s)=
R(1,\frac{s}{r})r^{2HK} $ and the change of variables $r/s= z$ in
the integral respect to $dz$ to obtain
\begin{eqnarray*}
\|J_i^{\varepsilon}(t)\|_{\alpha-1,2}^2&\leq&C\sum_{m\geq0}(1+m)^{\alpha}\sum_{|n|=n_1+...+n_d=m}\prod_{j=1}^{d}\frac{1}{(n_j\vee1)^{\frac{8\beta-1}{6}}}
\\&\times&
\int_0^tr^{2(\gamma-1+H_iK_i)}dr\int_0^1(z)^{\gamma-\frac{3}{2}+H_iK_i}\left(\frac{
R_{i}(1,z) }{z^{H_{i}K_{i}}}\right) ^{n_{i}-1}\prod_{\underset{j\neq
i}{j=1}}^{d} \left( \frac{ R_{j}(1,z) }{z^{H_{j}K_{j}}}\right)
^{n_{j}}dz.
\end{eqnarray*}
Since for each $i\in\{1,\ldots,d\}$, $\gamma-\frac{3}{2}+H_iK_i>-1$,
$2(\gamma-1+H_iK_i)>-1$ and from Lemma 12 and  the proof of the
Proposition 12 in \cite{ELSTV}, we obtain that
\begin{eqnarray*}
\|J_i^{\varepsilon}(t)\|_{\alpha-1
,2}^2&\leq&C\sum_{m\geq0}(1+m)^{\alpha} m^{-\frac{1}{2(HK)^{\ast} }}
\sum_{|n|=n_1+...+n_d=m}\prod_{j=1}^{d}\frac{1}{(n_j\vee1)^{\frac{8\beta-1}{6}}}\\
&\leq &C\sum_{m\geq0}(1+m)^{\alpha}
(m)^{-\frac{1}{2(HK)^{\ast}}-1+d(1-\frac{8\beta -1}{6})}
\end{eqnarray*}
and since is finite if and only if
$\alpha<\frac{1}{2(HK)^*}-\frac{d}{2}$.\\ On the other hand, for
every $(s,z)\in(0,T]\times\R^d$ we have
\begin{eqnarray*}\left|\partial_s
\bar{U}(s,z)\right|\leq C
s^{\frac{-d}{2}+\theta}\left|\left((z_1-x_1)s^{-H_1K1},...,(z_d-x_d)s^{-H_dKd}\right)\right|^{2-d}.
 \end{eqnarray*}and
 \begin{eqnarray*}\left|
\bar{U}(s,z)\right|\leq C
s^{\frac{1}{2}(2-d)+\theta}\left|\left((z_1-x_1)s^{-H_1K1},...,(z_d-x_d)s^{-H_dKd}\right)\right|^{2-d}
 \end{eqnarray*}
 and this inequalities imply as in \cite{uemura} the convergences
(\ref{converg.ofpartial_sUin D}) and (\ref{converg. of U in D}) in $mathbb{D}^{\alpha}_2$.\qed

\bibliographystyle{amsplain}

\end{document}